\documentclass{article}
\newcommand{\figd}{./png}
\usepackage{amsmath}
\usepackage{fancyhdr}
\usepackage{graphicx}

\newcommand{\thup}{\theta_{\uparrow}}
\newcommand{\thdn}{\theta_{\downarrow}}

\newcommand{\dthi}{\partial_{\theta_i}}
\newcommand{\dthj}{\partial_{\theta_j}}
\newcommand{\onehalf}{{\scriptstyle\frac{1}{2}}}
\newcommand{\ul}{\underline}
\renewcommand{\th}{^{\mbox{\ul{\scriptsize th}}}}

\title{The Optimal `AND'}
\author{Richard Rohwer}
\date{\today}

\begin{document}
\maketitle\thispagestyle{fancy}
\begin{abstract}
  The joint distribution $P(X,Y)$ cannot be determined from its
  marginals $P(X)$ and $P(Y)$ alone; one also needs one of the
  conditionals $P(X|Y)$ or $P(Y|X)$.  But is there a best guess,
  given only the marginals?  Here we answer this question in the
  affirmative, obtaining in closed form the function of the marginals
  that has the lowest expected 
  Kullbach-Liebler (KL) divergence between the unknown ``true''
  joint probability and the function value.
  The expectation is taken with respect to
  Jeffreys' non-informative prior over the possible joint
  probability values, given
  the marginals.  This distribution can also be used to obtain the
  expected information loss for any other {\em aggregation operator},
  as such estimators are often called in fuzzy logic, for any given
  pair of marginal input values.  This enables such such operators,
  including ours, to be compared according to their expected loss
  under the minimal knowledge conditions we assume.

  We go on to develop a method for evaluating the expected accuracy of
  any aggregation operator in the absence of knowledge of its inputs.
  This requires averaging the expected loss over all possible input
  pairs, weighted by an appropriate distribution.  We obtain this
  distribution by marginalizing Jeffreys' prior over the possible
  joint distributions (over the 3 functionally
  independent coordinates of the space of joint distributions over two Boolean
  variables) onto a joint distribution over the pair of marginal
  distributions, a 2-dimensional space with one parameter for each
  marginal.  We report the resulting input-averaged expected losses
  for a few commonly used operators, as well as the optimal
  operator.

  Finally, we discuss the potential to develop our methodology into
  a principled risk management approach to replace
  the often rather arbitrary conditional-independence assumptions
  made for probabilistic graphical models.
\end{abstract}

\section{Introduction}
Given truth values for propositions $X$ and $Y$, the truth value of
the conjunction
$X\wedge Y$ is fully determined by the truth table for the logical
connective `AND'.  But if we are only able to assign non-extremal
probabilities $P(X)$ and $P(Y)$ to these propositions, then the
joint probability $P(X\wedge Y)$ is not fully determined.
Further information such as the conditional probability $P(Y|X)$
is required.  The essence of the problem is that the space of joint
distributions over two Boolean values has three dimensions, of which 
the marginals are only able to constrain two.
Nevertheless, some guesses are better than others.
Under Bayesian principles~\cite{jaynes03,Caticha:2008eso},
the best guess is the one that minimizes the
expected information loss under the distribution that maximizes
the entropy over what is not known, given whatever is known.
Here we apply this principle
to obtain an optimal estimate for $P(X\wedge Y)$ given $P(X)$
and $P(Y)$.  We apply the principle again under the still more
information-sparse conditions of knowing neither $P(X)$ nor $P(Y)$
to obtain the expected information loss for any given aggregation
operator that estimates $P(X\wedge Y)$ from these marginals.

It is important to notice that we are not directly concerned with
whether propositions $X$ and $Y$ hold true; the incompletely known
variables in our chosen problems are probability distributions
over the four possible joint outcomes of these two Boolean variables.
To describe our state of knowledge about these distributions, we
must work at the meta-level in which events take values within the
space of such distributions, and concern ourselves with distributions
over these distribution spaces.  Unlike a finite set, which supports a
trivial concept of uniformity, this event space forms a continuum, so entropy
can be defined only with respect to some fiducial distribution that
{\em defines} uniformity.  In spaces of distributions, the uniquely qualified
candidate for this meta-distribution is Jeffreys' non-informative prior.
It is proportional to the square root of determinant of the Fisher
information metric, which in turn measures the distinguishibility
of infinitesimally nearby distributions in the space.
A crash course in these information-geometric topics can be
found in Appendix~\ref{sec:appA}.
Given this apparatus, maximum entropy problems become well-formulated
as minimum Kullbach-Liebler divergence (KL) problems.

\section{Formulation}
We begin by introducing coordinates $\theta$ for the space of joint
categorical distributions over two Boolean variables:
\begin{align}
  \theta_{11}&=P(X,Y)       \label{eq:theta11} \\
  \theta_{10}&=P(X,\neg Y)  \label{eq:theta10} \\
  \theta_{01}&=P(\neg X, Y) \label{eq:theta01}
\end{align}
The coordinates range over the 3-simplex defined by
\begin{align}
  \theta_{11} &\geq 0 \label{eq:theta11pos}\\
  \theta_{10} &\geq 0 \label{eq:theta10pos}\\
  \theta_{01} &\geq 0 \label{eq:theta01pos}\\
  \theta_{11} + \theta_{10} + \theta_{01} &\leq 1 \label{eq:thetaNorm}
\end{align}
We also introduce the function
\begin{align}
  \theta_{00}(\theta)=1-\theta_{11} - \theta_{10} - \theta_{01}
\end{align}
which is not to be confused with one of the 3 coordinates,
even when written as $\theta_{00}$ without making its argument explicit.

In this coordinate system, we obviously have
\begin{align}
  P(X) &= \theta_{11}+\theta_{10} \\
  P(Y) &= \theta_{11}+\theta_{01}
\end{align}

Abbreviating $P(X)$ as $a$ and $P(Y)$ as $b$, our knowledge can be expressed
by the constraints
\begin{align}
  \theta_{11}+\theta_{10} &= a \label{eq:aConstr}\\
  \theta_{11}+\theta_{01} &= b \label{eq:bConstr}
\end{align}

Let us use these constraints to eliminate the coordinates
$\theta_{10}$ and $\theta_{01}$:
\begin{align}
  \theta_{10} &= a - \theta_{11}\label{eq:aConstr1}\\
  \theta_{01} &= b - \theta_{11}\label{eq:bConstr1}\\
\end{align}

The range restrictions (\ref{eq:theta11pos}), (\ref{eq:theta10pos}),
(\ref{eq:theta01pos}) and (\ref{eq:thetaNorm}) then become
\begin{align}
  \theta_{11} &\geq 0 \label{eq:theta11posA} & \theta_{11} &\geq 0 \\
  a - \theta_{11} &\geq 0 \label{eq:theta10posA} &\qquad
  \theta_{11} &\leq a \\
  b - \theta_{11} &\geq 0 \label{eq:theta01posA} &\qquad
  \theta_{11} &\leq b &\qquad  \\
  \theta_{11}+a-\theta_{11}+b-\theta_{11} &\leq 1 &\qquad
  \theta_{11} &\geq a+b-1   \label{eq:thetaNormA}
\end{align}
which can be summarized as
\begin{align}
  \thdn & \leq \theta_{11} \leq \thup \label{eq:Th11range}\\
  \thdn &= \max(a+b-1,0) \label{eq:ThUp}\\
  \thup &= \min(a,b) \label{eq:ThDn}\\
\end{align}

\section{The optimal estimate}
\begin{figure}
  \begin{center}
    \begin{tabular}{cc}
      \includegraphics[scale=0.4,trim=0 0 30 0,clip]{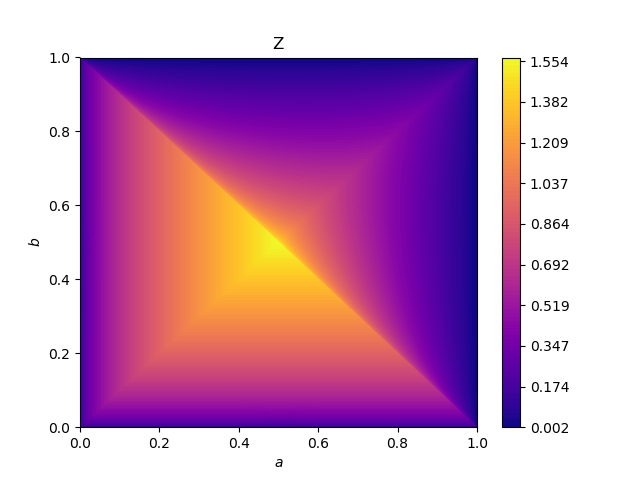} &
      \includegraphics[scale=0.5,trim=80 20 0 0,clip]
                      {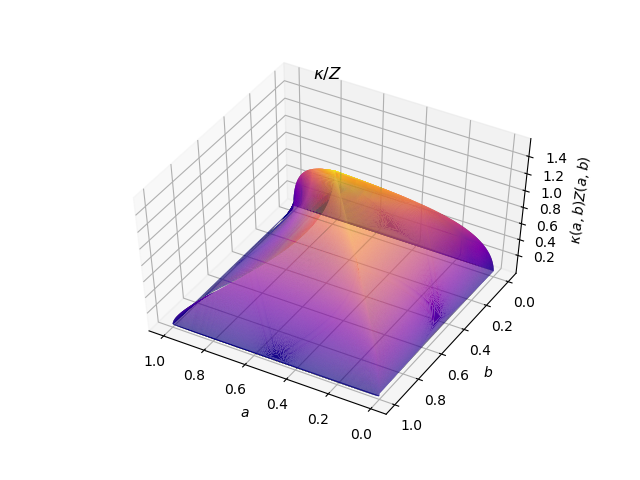}
    \end{tabular}
  \end{center}
  \caption{The normalization factor $Z$ (\ref{eq:Zdef}) as a function of
    the marginals.  The discontinuities of (\ref{eq:ThUp})
    and (\ref{eq:ThDn}) are prominently reflected.
  }
\label{fig:Z}
\end{figure}

In Appendix~\ref{sec:appA}, we motivate and derive the Fisher information
matrix and the measure that it defines over spaces of probability
distributions.  In the present case,
the Fisher information matrix has the single element
\begin{align}
  g_{(11),(11)}(\theta_{11}) &= 
  \sum_{x\in\{T,F\}} P(x|\theta_{11})
  [\partial_{\theta_{11}} \ln P(x|\theta_{11})]
  [\partial_{\theta_{11}} \ln P(x|\theta_{11})] \notag\\
  &= \theta_{11} \theta_{11}^{-2}
  + (1-\theta_{11}) (-(1-\theta_{11}))^{-2}\notag\\
  &= \theta_{11}^{-1} + (1-\theta_{11})^{-1}
  = 1/(\theta_{11} (1-\theta_{11}))
\end{align}
The ``uniform'' measure over this space, Jeffreys' non-informative prior,
is the square root of the determinant of $g$, which is just
$1/\sqrt{\theta_{11} (1-\theta_{11})}$.  This is essentially the
Beta($\onehalf,\onehalf$) density, except that its support is limited
to a sub-interval of the unit interval.  The integral over this
support region, which gives the normalization constant needed to
convert this measure into a probability measure, is 
an incomplete Beta function for which we have a special case that
reduces to
\begin{align}
  Z &= \int_{\thdn}^{\thup}\frac{dx}{\sqrt{x (1-x)}} =
  2 \sin^{-1}(\sqrt{x})\biggr|_{\thdn}^{\thup} \notag\\
  &= 2 \sin^{-1}\bigl(\sqrt{\thup}\bigr)
  - 2\sin^{-1}\bigl(\sqrt{\thdn}\bigr) \label{eq:Zdef}
\end{align}
which is plotted in Figure~\ref{fig:Z}.

To see this, recall the standard result
$d\sin^{-1}(x)/dx=1/\sqrt{1-x^2}$ and observe that
\begin{align}
  \frac{d}{dx} 2 \sin^{-1}(\sqrt{x}) =
  \frac {2}{\sqrt{1-x}}\frac{1}{2\sqrt{x}}
  = \frac{1}{\sqrt{x(1-x)}}.
\end{align}
We therefore conclude that $\theta_{11}$ is distributed between
$\thdn$ and $\thup$ according to the density
\begin{align}
  P(\theta_{11}) = \frac{1}{Z}
  \frac{1}{\sqrt{(\theta_{11} (1-\theta_{11}))}}  \label{eq:P11}
\end{align}
with $Z$ given by (\ref{eq:Zdef}).

With this distribution in hand, we can derive an optimal estimate
$\theta_{11}^*(a,b)$ for $\theta_{11}$.  Suppose we institute a policy
of using some fixed function $f(a,b)$ to estimate $\theta_{11}$.
Then the expected KL divergence between this estimate and
$\theta_{11}$ would be
\begin{align}
  KL[\theta_{11},f(a,b)] &=
  \theta_{11}\ln\frac{\theta_{11}}{f(a,b)} +
  (1-\theta_{11})\ln\frac{1-\theta_{11}}{1-f(a,b)} \label{eq:KL11}
\end{align}
where we choose to weight the log likelihood ratio
by the true distribution $\theta_{11}$, not
the estimated distribution $f(a,b)$.
The expected information loss is therefore the expectation value of 
(\ref{eq:KL11}) under distribution (\ref{eq:P11}).
Abbreviating $f(a,b)$ as $x$ and $\theta_{11}$ as $\xi$,
this is
\begin{align}
  \langle KL[\theta_{11},f(a,b)]\rangle_{\theta_{11}} &=
  \frac{\kappa}{Z} - \frac{A}{Z}\ln x - \frac{B}{Z}\ln(1-x) \label{eq:ExpKL}\\
  \kappa &= \int_{\thdn}^{\thup} d\xi
  \frac{\xi\ln\xi + (1-\xi)\ln(1-\xi)}{\sqrt{\xi(1-\xi)}} \label{eq:Kappadef}\\
  A &= \int_{\thdn}^{\thup} d\xi \sqrt{\xi/(1-\xi)} \label{eq:Adef}\\
  B &= \int_{\thdn}^{\thup} d\xi \sqrt{(1-\xi)/\xi} \label{eq:Bdef}\\
\end{align}
Note that with the change of variable $\xi'=1-\xi$, $d\xi'=-d\xi$,
(\ref{eq:Bdef}) becomes
$B=-\int_{1-\thdn}^{1-\thup} d\xi'\sqrt{\xi'/(1-\xi')}=
\int_{1-\thup}^{1-\thdn} d\xi'\sqrt{\xi'/(1-\xi')}$, so $A$ and $B$
involve the same indefinite integral, verified in Appendix~\ref{sec:intglA}
(\ref{eq:intglAchk}) to be:
\begin{align}
  \int d\xi \sqrt{\frac{\xi}{1-\xi}} = \tan^{-1}(\sqrt{\xi/(1-\xi)})
  -\sqrt{\xi(1-\xi)}. \label{eq:intglA}
\end{align}
This gives
\begin{align}
  A =& \tan^{-1}\Bigl(\sqrt{\thup/(1-\thup)}\Bigr)
  - \tan^{-1}\Bigl(\sqrt{\thdn/(1-\thdn)}\Bigr) \notag\\ &
  -\sqrt{\thup(1-\thup)} +\sqrt{\thdn(1-\thdn)}.  \label{eq:A}\\
  B =& \tan^{-1}\Bigl(\sqrt{(1-\thdn)/\thdn}\Bigr)
  - \tan^{-1}\Bigl(\sqrt{(1-\thup)/\thup}\Bigr) \notag\\ &
  -\sqrt{(1-\thdn)\thdn} +\sqrt{(1-\thup)\thup}.  \label{eq:B}
\end{align}
We note that $A$ and $B$ are both non-negative because the integrand
of (\ref{eq:intglA}) is non-negative and in each case the integral
is taken in a positive sense.

\begin{figure}
  \begin{center}
    \begin{tabular}{cc}
      \includegraphics[scale=0.4,trim=0 0 30 0,clip]{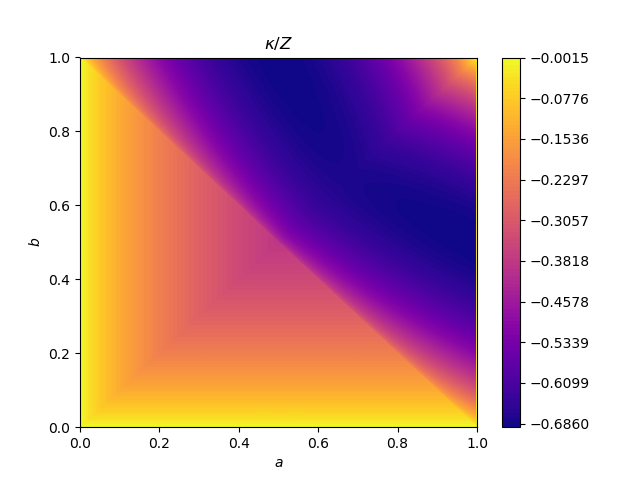} &
      \includegraphics[scale=0.5,trim=80 20 0 0,clip]
                      {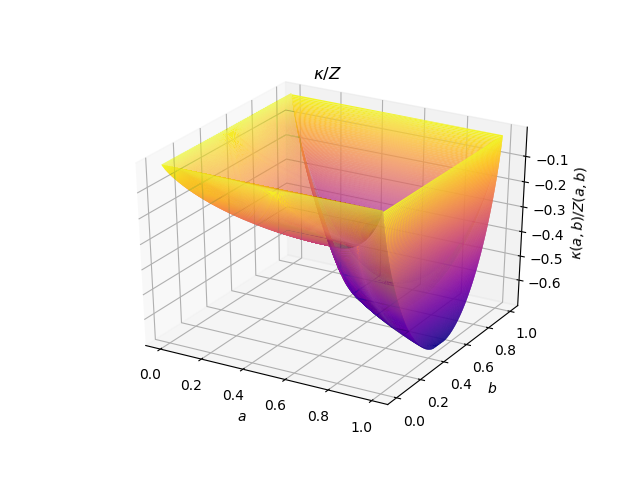}
    \end{tabular}
  \end{center}
  \caption{The aggregation-operator-independent term $\kappa/Z$
    (\ref{eq:Kappadef},\ref{eq:Zdef}) as a function of the input marginals.
    Contour plot on the right, 3d plot on the left.
  }
\label{fig:kappa}
\end{figure}

The two terms in $\kappa$ (\ref{eq:Kappadef}) can be similarly written
in terms of a single indefinite integral, but one that does not have
a simple expression in terms of elementary functions.  (It can be
expressed in terms of a generalized hypergeometric function, but
we do not pursue this here.)
It was integrated numerically to obtain the plot
in Figure~\ref{fig:kappa}.  We note that $\kappa\leq 0$, as is clear
from (\ref{eq:Kappadef}).
This term does not depend on the
estimate $x=f(a,b)$, and in this sense is a constant
contribution to the expected information loss.

The optimal estimate $\theta_{11}^*(a,b)$ minimizes the expected loss
(\ref{eq:ExpKL}).
Setting its derivative with respect to $x$ to zero gives
\begin{align}
  \frac{A}{x} - \frac{B}{1-x} &= 0 \notag\\
  A - Ax &= Bx \notag\\
  x &= \theta_{11}^*(a,b) = A/(A+B) \label{eq:optTh}
\end{align}
We note that $A+B$ contains only the $\tan^{-1}$ terms of
(\ref{eq:A}) and (\ref{eq:B}).

Because $A$ and $B$ are non-negative, we see that the optimal estimate
(\ref{eq:optTh}) lies between 0 and 1, as it should.  To confirm
that the extremum (\ref{eq:optTh}) is indeed a minimum, we note 
that the second derivative of (\ref{eq:ExpKL}) is proportional to
\begin{align}
  \frac{A}{x^2} + \frac{B}{(1-x)^2} \geq 0.
\end{align}

\begin{figure}
  \begin{center}
    \begin{tabular}{cc}
      \includegraphics[scale=0.4,
        trim=0 0 50 0,clip]{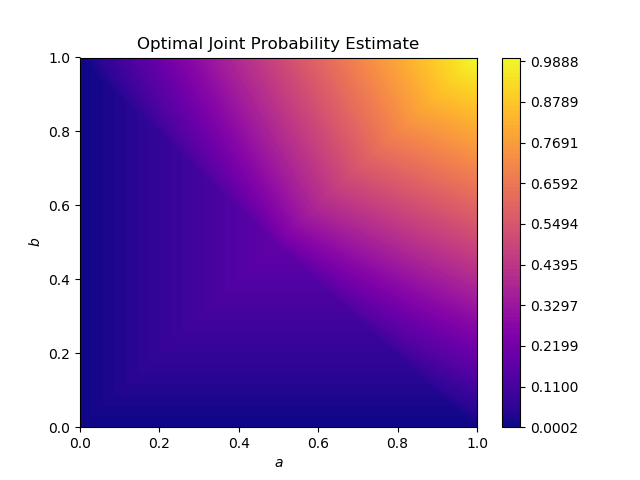} &
      \includegraphics[scale=0.5,
        trim=90 20 0 60,clip]{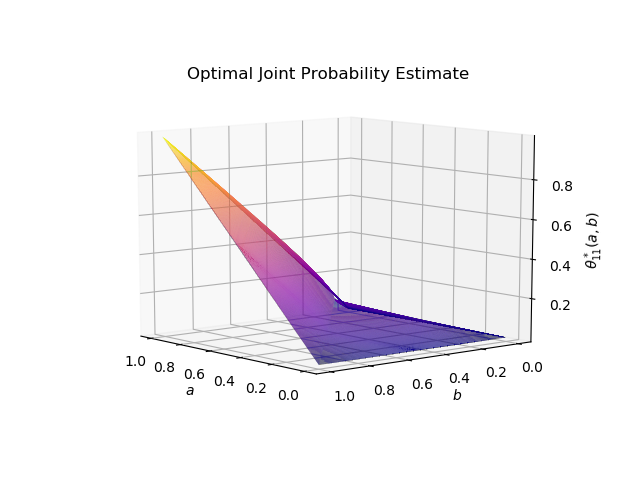} \\
      \includegraphics[scale=0.4,
        trim=0 0 50 0,clip]{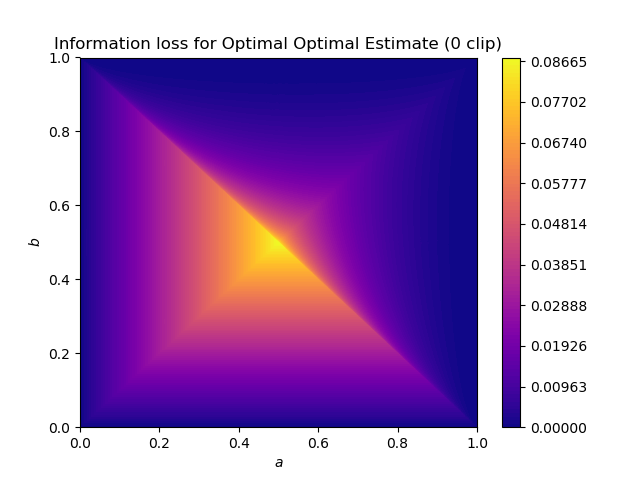} &
      \includegraphics[scale=0.5,
        trim=90 20 0 60,clip]{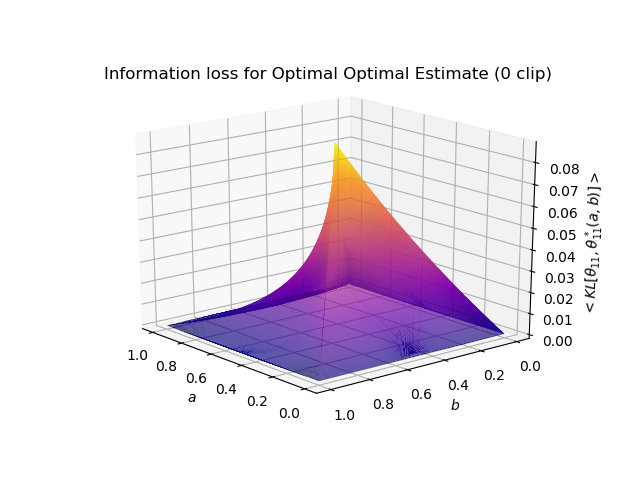}
    \end{tabular}
  \end{center}
  \caption{The optimal aggregation operator (\ref{eq:optTh}).
    Contour plot on upper left, 3D
    plot on upper right.  Corresponding information loss costs in lower
    two plots.  These include the aggregation-operator-independent
    term $\kappa$ (\ref{eq:Kappadef}) shown in Figure~\ref{fig:kappa},
    which was obtained by numerical integration.  Small inaccuracies
    in this integration resulted in small negative values for the
    expected KL cost, which are not mathematically possible.  These
    were truncated to 0.}
\label{fig:thetaOptCopt}
\end{figure}

Plugging $x$ from (\ref{eq:optTh}) into (\ref{eq:ExpKL}) gives the
expected loss for the optimal estimate:
\begin{align}
  \langle KL[\theta_{11},\theta_{11}^*(a,b)]\rangle_{\theta_{11}} &=
  \frac{\kappa}{Z} - \frac{A}{Z}\ln \frac{A}{A+B}
  - \frac{B}{Z}\ln\frac{B}{A+B} \notag\\
  &= \bigl(\kappa -A\ln A -B\ln B +(A+B)\ln(A+B)\bigr)/Z
  \label{eq:ExpKLopt}
\end{align}
The optimal estimate (\ref{eq:optTh}) is plotted along with its expected
loss (\ref{eq:ExpKLopt}) in Figure~\ref{fig:thetaOptCopt}.

\section{Aggregation operator cost}
Expression (\ref{eq:optTh}), together with (\ref{eq:A}), (\ref{eq:B}),
(\ref{eq:ThUp}) and (\ref{eq:ThDn}) gives the optimal estimate
$P(X,Y)$ given that $P(X)=a$ and $P(X)=b$.
Given a different estimate $f(a,b)$, of which a considerable variety
are in use~\cite{Detyniecki:01:Aggregation},
one can apply (\ref{eq:ExpKL}) and compare the
result to (\ref{eq:ExpKLopt}) to determine how much worse the
sub-optimal estimate is for any particular marginals $a$ and $b$.

In order to judge how well an aggregation formula $f$ works in general,
without committing to any particular values of its arguments,
we need to average the expected loss of the estimate
over all marginal values $a$ and $b$ under a non-informative
distribution.  To obtain this distribution over $(a,b)$, we
assume Jeffreys' non-informative prior over the joint distributions
$\theta=(\theta_{11},\theta_{01},\theta_{10})$, apply 
(\ref{eq:aConstr}) and (\ref{eq:bConstr}) to transform to
coordinates that eliminate $\theta_{01}$ and $\theta_{10}$
in favor of $a$ and $b$, and integrate out $\theta_{11}$.

As derived in Appendix~\ref{sec:appA}, Jeffreys' prior over the joint
distributions is the Dirichlet distribution:
\begin{align}
  P(\theta) =
  \pi^{-2} (\theta_{11}\theta_{10}\theta_{01}\theta_{00})^{-\onehalf}
  \label{eq:Jprior3}
\end{align}

To change to coordinates $\phi\equiv(\theta_{11},a,b)$, we must not
only substitute (\ref{eq:aConstr1}) and (\ref{eq:bConstr1}) into
(\ref{eq:Jprior3}), but also ensure that
$P(\theta)d^3\theta=P(\phi)d^3\phi$ by dividing $P(\theta)$ by the 
absolute value of the Jacobian determinant
\begin{align}
  \begin{vmatrix}
    \frac{\partial\phi}{\partial\theta}
  \end{vmatrix}
   =
  \begin{vmatrix}
  \frac{\partial\theta_{11}}{\partial\theta_{11}} &
  \frac{\partial a}{\partial \theta_{11}} & 
  \frac{\partial b}{\partial \theta_{11}} \\
  \frac{\partial\theta_{11}}{\partial \theta_{01}} &
  \frac{\partial a}{\partial \theta_{01}} & 
  \frac{\partial b}{\partial \theta_{01}} \\
  \frac{\partial\theta_{11}}{\partial \theta_{10}} &
  \frac{\partial a}{\partial \theta_{10}} & 
  \frac{\partial b}{\partial \theta_{10}}
  \end{vmatrix}
  =
  \begin{vmatrix}
    1 & 1 & 1 \\
    0 & 0 & 1 \\
    0 & 1 & 0
  \end{vmatrix}
  =-1.
\end{align}
However, this turns out to be an immaterial unit factor, so we have
\begin{align}
  P(\phi) &= \pi^{-2}\bigl(\theta_{11}(a-\theta_{11})(b-\theta_{11})
  (1-\theta_{11}-(a-\theta_{11})-(b-\theta_{11})\bigr)^{-\onehalf} \notag\\
  &= \pi^{-2}\bigl(\theta_{11}(a-\theta_{11})(b-\theta_{11})
  (1+\theta_{11}-a-b)\bigr)^{-\onehalf}
\end{align}
The marginal over $a$ and $b$ is then found by integrating out $\theta_{11}$:
\begin{align}
  P(a,b) &= \pi^{-2} \int_{\thdn}^{\thup} d\theta_{11}
   \bigl(\theta_{11}(a-\theta_{11})(b-\theta_{11})
  (1+\theta_{11}-a-b)\bigr)^{-\onehalf}. \label{eq:abMrg}
\end{align}
This is an incomplete elliptic integral of the first kind.
With the definition
\begin{align}
  K(\phi,m) = \int_0^{\phi} (1-m\sin^2\theta)^{-\onehalf} d\theta
  \qquad 0\leq m \leq 1,\quad 0\leq\phi\leq\frac{\pi}{2}
  \label{eq:elipDef}
\end{align}
and substituting $x=\theta_{11}$,
the antiderivative of the integrand 
$\pi^{-2}\bigl(x(a-x)(b-x)(1+x-a-b)\bigr)^{-\onehalf}$
in (\ref{eq:abMrg}) is 
\begin{align}
  I(x;a,b) = \frac{2}{\pi^2}\frac{b-a}{|b-a|} \frac{1}{\sqrt{b(1-b)}}
  K \Biggl(\sin^{-1}\Biggl(\sqrt{\frac{(1-b)(b-x)}{(1-a)(a-x)}}\Biggr),
  \frac{a(1-a)}{b(1-b)}\Biggr) \label{eq:indefIntEllip}
\end{align}
as verified in Appendix~\ref{sec:indefIntEllip}.
However, if we attempt to obtain
the definite integral (\ref{eq:abMrg}) simply by plugging in the limits,
we may find ourselves in violation of the conditions $0\leq m<1$
and $0\leq\phi\leq\pi/2$.  Inspecting (\ref{eq:indefIntEllip}), we see that
\begin{align}
  m &= \frac{a(1-a)}{b(1-b)} \\
  \sin^2\phi &= \frac{(1-b)(b-x)}{(1-a)(a-x)}.
\end{align}
The non-negativity constraint on $m$ is always satisfied, because $a$ and
$b$ are probabilities lying between $0$ and $1$.  But the upper bound is
problematic when $b$ is near either end of its allowed range.  That condition
can be written
\begin{align}
  \frac{a(1-a)}{b(1-b)} &< 1 \notag\\
  a-a^2 &< b-b^2 \notag\\
  0 &< b-a -(b^2-a^2) = (b-a)-(b-a)(b+a) \notag\\
  0 &< (b-a)(1-b-a) \notag\\
  (b>a \text{ and } a+b<1) &\text{ or } (b<a \text{ and } a+b>1)
  \label{eq:mConditions}
\end{align}
The remaining condition amounts to 
$0\leq\frac{(1-b)(b-x)}{(1-a)(a-x)}\leq 1$.  Because $a$ and $b$ are
marginal probabilities that include $x$ ({\em i.e.}, $\theta_{11}$)
within their mass, we have $b-x\geq 0$ and $a-x\geq 0$, so
the non-negativity constraint is satisfied.  The remaining constraint is
\begin{align}
  \frac{(1-b)(b-x)}{(1-a)(a-x)} &\leq 1 \notag\\
  a-x-a^2+ax &\geq b-x-b^2+bx \notag\\
  0 &\geq (b-a) - (b^2-a^2) +(b-a)x = (b-a)(1-a-b+x) \notag\\
  0 & \leq (b-a)(a+b-1-x) \notag\\
  (b\geq a \text{ and } x \leq a+b-1) &\text{ or }
  (b \leq a \text{ and } x\geq a+b-1)
  \label{eq:phiConditions}
\end{align}
Combining these conditions with (\ref{eq:mConditions}) has the consequences
\begin{align}
  (b\geq a \text{ and } x\leq 0) \text{ or } (b\leq a \text{ and } x\geq 0)
\end{align}
so because $x\geq 0$ we conclude that we can proceed only if $b\leq a$,
in which case (\ref{eq:mConditions}) implies that we must also
restrict consideration to $a+b>1$.

It remains to determine whether the second condition in
(\ref{eq:phiConditions}), $x\geq a+b-1$, is respected by the
limits $x=\thup=\min(a,b)$ and $x=\thdn=\max(a+b-1,0)$.
For $\thup$, having already restricted to $b\leq a$, 
we have
\begin{align}
  \min(a,b) &\geq a+b-1 \notag\\
  b &\geq a+b-1 \notag\\
  1 &\geq a
\end{align}
which is always satisfied.  Moving on to $\thdn$ we obtain
\begin{align}
  \max(a+b-1,0) &\geq a+b-1.
\end{align}
which, having restricted attention to $a+b>1$, is always satisfied.

Leaving aside the question of how to handle the remaining cases,
let us plug the limits into (\ref{eq:indefIntEllip}).
At $x=\thup=\min(a,b)=b$, the
argument of $K$ vanishes, so as is clear from (\ref{eq:elipDef}),
so does $K$.  This leaves the lower limit
$x=\thdn=\max(a+b-1,0)=a+b-1$ (with the $a+b>1$ constraint), at
which we see
\begin{align}
  \frac{(1-b)(b-x)}{(1-a)(a-x)} &=
  \frac{(1-b)(b-(a+b-1))}{(1-a)(a-(a+b-1))} =
  \frac{(1-b)(1-a)}{(1-a)(1-b)} = 1
\end{align}
Thus, we obtain
\begin{align}
  P(a,b) = \frac{2}{\pi^2} \frac{1}{b(1-b)}
  K\Bigl(\frac{\pi}{2},\frac{a(1-a)}{b(1-b)}\Bigr) \label{eq:abPrior0}
\end{align}
for $b<a$ and $a+b>1$.  Let us call this expression $P^{<>}(a,b)$
to remind us of these two conditions.

The remaining cases can be obtained from the symmetries of 
(\ref{eq:abMrg}).  The most obvious is that $P(a,b)=P(b,a)$.
Hence, for all $a$ and $b$ with $a+b>1$, we can say
\begin{align}
  P(a,b)=
  \begin{cases}
    P^{<>}(a,b) & [b<a] \\
    P^{<>}(b,a) & [b>a]
  \end{cases} \label{eq:abPrior}
\end{align}
leaving the $a+b<1$ cases to be determined.

\begin{figure}
  \begin{center}
    \begin{tabular}{cc}
      \includegraphics[scale=0.4,
        trim=0 0 50 0,clip]{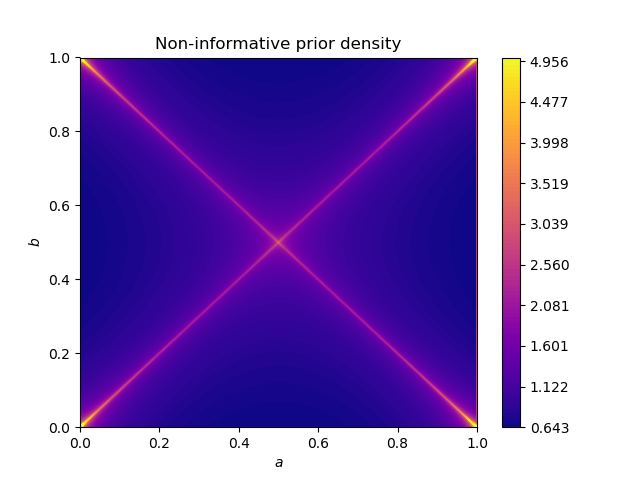} &
      \includegraphics[scale=0.5,
        trim=90 20 0 60,clip]{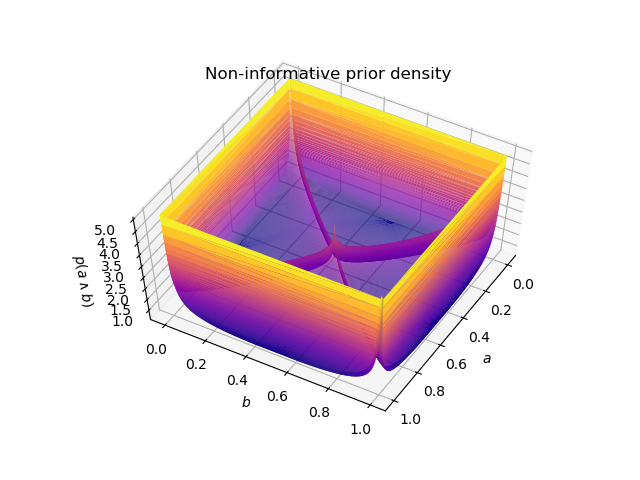}
    \end{tabular}
  \end{center}
  \caption{The non-informative prior distribution over the
    marginals. (\ref{eq:reflectSym}).
    Contour plot the left, 3D plot the right.  The distribution
    is singular at every edge, and along the diagonal lines.}
\label{fig:Pab}
\end{figure}

Less obviously, $P(1-a,1-b)=P(a,b)$.  Plugging into (\ref{eq:abMrg}),
we have
\begin{align}
  P(1-a,1-b) &= \pi^{-2} \int_{\thdn}^{\thup} d\theta_{11}
  \bigl(\theta_{11}(1-a-\theta_{11})(1-b-\theta_{11}) \notag\\
   &\qquad\qquad\qquad\qquad
  (1+\theta_{11}-(1-a)-(1-b))\bigr)^{-\onehalf} \\
  &= \pi^{-2} \int_{\thdn}^{\thup} d\theta_{11}
  \bigl(\theta_{11}(1-a-\theta_{11})(1-b-\theta_{11})
  (a+b+\theta_{11}-1)\bigr)^{-\onehalf} \notag
\end{align}
Let $x'=x+a+b-1$, so $x=x'-a-b+1$.  At the limit $x=\thup=\min(1-a,1-b)$
we have $x'=\min(1-a,1-b)+a+b-1$.  If $a<b$, this gives $x'=a$, and if
$b<a$ it gives $x'=b$, so this limit becomes $\min(a,b)$, as in
(\ref{eq:abMrg}).  At the limit
$x=\thdn=\max((1-a)+(1-b)-1,0)=\max(1-a-b,0)$, we have
$x'=\max(1-a-b,0)+a+b-1$.  If $a+b>1$, this is $a+b-1$, and if
$a+b<1$ it is $0$, so this limit becomes $\max{a+b-1,0}$, also
agreeing with (\ref{eq:abMrg}).  Thus
\begin{align}
  P(1-a,1-b) &= \pi^{-2} \int_{\thdn}^{\thup} dx
  \bigl(x(1-a-x)(1-b-x) (a+b+x-1)\bigr)^{-\onehalf} \notag\\
  &=
  \pi^{-2} \int_{\thdn}^{\thup} dx'
  \bigl( (x'-a-b+1) (1-a-(x'-a-b+1)) \notag\\
  & \qquad\qquad
  (1-b-(x'-a-b+1)) (a+b+(x'-a-b+1)-1)\bigr)^{-\onehalf} \notag\\
  &=
  \pi^{-2} \int_{\thdn}^{\thup} dx' (1+x'-a-b)(b-x')(a-x')x' \notag\\
  &=P(a,b)  \label{eq:reflectSym}
\end{align}
If $a+b<1$, then $(1-a)+(1-b)=2-(a+b)>1$, so we can use
(\ref{eq:reflectSym}) to obtain $P(a,b)$.
This density is plotted in Figure~\ref{fig:Pab}, which clearly shows
the singularities at $a=b$ and $a+b=1$.

\begin{table}[!h]
  \begin{center}
  \begin{tabular}{|lll|}
    \hline
    {\bf Method} & {\bf Formula} & {\bf Expected loss} \\
    \hline
    Optimal & Eqn. (\ref{eq:optTh}) & 0.0203 \\
    Product & $ab$ & 0.0272 \\
    Hamacher ($\gamma=0$) & $\frac{ab}{a+b-ab}$ & 0.0416 \\
    Min & $\min(a,b)$ & 0.0691 \\
    Lukasiewicz &  $\max(a+b-1,0)$ & 0.5183 \\
    \hline
  \end{tabular}
  \caption{Expected information loss for the optimal estimate and a few
  widely used t-norm aggregation operators.}
  \label{tbl:league}
  \end{center}
\end{table}

The expected information loss (\ref{eq:ExpKL}) was integrated over this
density using the optimal estimate (\ref{eq:optTh})
and a few other widely used t-norm
aggregation operators taken from \cite{Detyniecki:01:Aggregation}
to write $x$ as a function of $a$ and $b$.  The resulting overall
expected losses are shown in Table~\ref{tbl:league}.

Details of the numerical integration are given in Appendix~\ref{app:Num}.
The numerical procedure applied to the density
(\ref{eq:abPrior0}, \ref{eq:abPrior}) alone integrated to 0.9989,
so these figures should be accurate to roughly 0.1\%.

\section{Summary and future work}
We have set out and demonstrated a principled method for deriving
optimal aggregation formulas.  The best-guess full joint distribution
over two Boolean variables, given the marginal distributions, is
given in closed form by (\ref{eq:optTh}), and its expected information
loss, measured by Kullbach-Liebler divergence, is given by
(\ref{eq:ExpKLopt}).
Of course, if one has further information then the maximum entropy
problem should be formulated differently, resulting in an
improved guess.

We also set out and demonstrated a principled method for evaluating
aggregation operators by their expected information loss, averaging
over all their possible inputs.  We applied the method to produce
Table~\ref{tbl:league}, listing the expected loss for the
optimal aggregation operator and a few others.  The method can be
applied to extend the table to anyone's favorite aggregation operator.

In probabilistic graphical model applications~\cite{bishopbook06},
conditional independence assumptions
are used to estimate joint distributions.
These assumptions are often motivated
more by practical limitations than by plausibility.  
At least in principle, it should be a better bet to take a 
quantitative risk management approach to dealing with such limitations,
which is essentially what we have demonstrated in the simplest non-trivial
case.  One places the bet that has the lowest risk of being badly wrong,
in the sense of expected information loss.
Therefore it is of considerable interest to attempt to generalize
the method in at least two ways: (i) from joint distributions over
Boolean variables to joint-distributions over $n$-valued categorical
variables for arbitrary $n>2$; and (ii) from joint distributions over
two variables to distributions over greater numbers.  On the latter point,
with the machinery given here it is possible to treat the case of
$n$ Boolean variables by repeated application of the 2-variable formula
(\ref{eq:optTh}), but we see no good reason to believe that this
would give the same result as minimization of expected KL loss within
the space of $n$-variable distributions.  It would also be of interest
to extend the approach to commonly-used families of distributions over
continuous spaces.

% Another extension would be to generalize from joint distributions over
% Boolean variables to joint-distributions over $n$-valued categorical
% variables, for arbitrary $n>2$.
%  This would have an immediate application
% to differentiable co-clustering.  In this algorithm, a learnable
% function maps a data object $x$ to a categorical distribution
% over $n$ feature values, say $p_i(d)$, $0\leq i\leq n-1$.
% A data set organizes data objects into $U$ pairs, the $u\th$ of which
% we may call $(d_0^u,d_1^u)$.  Then a joint distribution
% $p_{ij}=U^{-1}\sum_up_i(d_0^u)p_j(d_1^u)$ is defined, and its mutual
% information is maximized with respect to the learnable parameters
% of the feature generators $p_i$.  Each summand, $p_i(d_0^u)p_j(d_1^u)$
% is a maximum-entropy estimate of the joint distribution given
% the marginals for sample $u$, but there is no good reason to
% choose the maximum-entropy distribution; one simply does not know
% what the joint distribution should be.  An optimal estimate based
% on Jeffreys' prior would be a better guess.  A severe anomoly arising
% from this incorrect use of a maximum-entropy assumption
% can be found by examining
% the degenerate case $d_0^u=d_1^u$ (for all $u$), which yields no mutual
% information at all but would be more reasonably expected to yield
% mutual information equal to the marginal entropy.

\section{Acknowledgments}
This work was carried out entirely with the author's own time and resources,
but benefitted from useful conversations with Elizabeth Rohwer and
his SRI colleagues John Byrnes, Andrew Silberfarb and others on
the Deep Adaptive Semantic Logic (DASL) team, supported in part by
DARPA contracts HR001118C0023 and HR001119C0108 and SRI internal
funding.  All views expressed herein are the author's own and are not
necessarily shared by SRI or the US Government.

\appendix
\section*{Appendices}

\section{Fisher Information and Jeffreys' prior} \label{sec:appA}
In this brief introduction to these topics in information
geometry~\cite{amari93,Caticha_2015},
we begin with the principle of {\em invariance}, which holds that
any formula for comparing probability distributions should produce
the same result under any one-to-one smooth invertable change of
random variables.  Otherwise the comparison would depend not purely
on the distributions, but also on the coordinates used to express
the random variables.  This leads to the $\delta$-divergences, among
which is the Kullbach-Leibler (KL) divergence.  Departing momentarily
from this thread, we then introduce the Fisher information matrix
and its determinant, and explain the sense in which this determinant measures
the number of distinguishable distributions in a given infinitesimal
coordinate volume element.  Normalized, this density of distributions
defines Jeffreys' prior, expressing that a ``randomly chosen'' distribution
is more likely to be chosen from a (coordinate) region dense in distinguishable
distributions than a region containing few such.
We then derive the Fisher information
matrix from the delta divergence between infinitesimally nearby distributions.
Finally, we derive these quantities for the categorical distributions
of interest here.

\subsection{Invariant divergences}
The condition that a functional $D[P,Q]$ of a pair of
distributions $P(X)$ and $Q(X)$ be invariant with respect to any
one-to-one smooth change of variables $X{\rightarrow}Y(X)$
is rather restrictive.  It implies that $D$ must be one of the 
$\delta$-divergences
\begin{align}
				D_{\delta}[P,Q] = \frac{1}{\delta(1-\delta)}
				\left[1-\int_x P(X=x)^{\delta} Q(X=x)^{1-\delta}\right]
				\label{eq:Def:deltaDiv}
\end{align}
or a function constructed from these divergences.  It is easy to see
that $D_{\delta}$ is invariant.  Recall that probability
densities\footnote{
We will often abbreviate $P(X=x)$ as $P(x)$ when this notational abuse
introduces little risk of confusion.  For that matter, we will often
be cavalier about notationally distinguishing between 
probabilities and densities.
}
transform as $P(y)=|{\partial}y/{\partial}x|^{-1}P(x)$
in order to preserve the probabilities
$P(y)dy=\{|{\partial}y/{\partial}x|^{-1}P(x)\}\{|{\partial}y/{\partial}x|dx\}=P(x)dx$.
It is much more complicated~\cite{amari93}
to prove that all invariant divergences are based on these forms.

For discrete random variables, the invariance is with respect to
uninformative subdivision of events.  That is, if
discrete outcome $x$ is replaced by two possible outcomes $x_1$ and
$x_2$, with $P(x_1)+P(x_2)=P(x)$, $Q(x_1)+Q(x_2)=Q(x)$ and
$P(x_1)/P(x_2)=Q(x_1)/Q(x_2)$, then $D_{\delta}[P,Q]$ is unchanged.
Carrying out a continuum transformation in
 a discretized
approximation also leads to this statement.

The KL and reverse-KL divergences are obtained in the 
$\delta\rightarrow{1}$ and $\delta\rightarrow{0}$ limits, due to the
identity\footnote{
$\partial_{\delta} x^{\delta} = \partial_{\delta} e^{\delta\ln{x}} = 
 x^{\delta}\ln{x} \rightarrow \ln{x}$
}
$\frac{x^{\delta}-1}{\delta}\rightarrow_{\delta\rightarrow{0}}\log{x}$
The symmetric $\delta=\onehalf$ case is the Hellinger distance.

\subsection{The Fisher Metric}
Infinitesimally, all the $\delta$-divergences reduce to the
Fisher metric (up to a constant scale factor).

For any family of probability distributions indexed by $n$-dimensional
parameter $\theta$, the distance between distributions $\theta$ and
$\theta+d\theta$ under the Fisher metric is 
\begin{align}
dD = \frac{1}{2}\sum_{ij}d\theta_ig_{ij}d\theta_j \label{eq:FisherDist}
\end{align}
 where $g$ is the Fisher information matrix
\begin{equation}
	g_{ij}(\theta) = \int_x P(x|\theta) \left[\dthi\ln P(x|\theta)\right]
	                            \left[\dthj\ln P(x|\theta)\right]
				\label{eq:Def:gij}
\end{equation}
It is straightforward to verify (\ref{eq:FisherDist})
and (\ref{eq:Def:gij}) by expanding 
(\ref{eq:Def:deltaDiv}) to second order in $d\theta$, and using the
conditions $\int_x{\dthi}P(x|\theta)=0$ which follow from the
normalization constraint $\int_xP(x|\theta)=1$.  This calculation is
carried out in section \ref{sec:FisherFromD}.

The Fisher distance $dD$ has a direct interpretation in terms of the
amount of IID data required to distinguish $P$ from $Q=P+dP$.  To see
this when the values of $x$ are discrete,
note that the typical log likelihood ratio between the probability of
$T$ samples from $P$, according to $Q$ and according to $P$ is
\begin{multline}
-\ln \left[ \prod_x Q(x)^{TP(x)} \right/\left. 
		 \prod_x P(x)^{TP(x)}\right] = \notag\\
		 T\sum_x P(x)\ln P(x)/Q(x) = TD_1[P,Q] 
\xrightarrow{Q\rightarrow P+dP} T dD.
\end{multline}
If we required this log likelihood ratio to
exceed some threshold $\alpha$ to declare $P$ and $P+dP$ distinct,
this condition would be $TdD>\alpha$ or $T>\alpha/dD$, so
up to a proportionality constant, $1/dD$ is the amount of data
required to make this distinction.  The result carries over to
any distribution over the continuum that is sufficiently regular to be
approximated by quantizing the continuum into cells.

It is shown directly in Section~\ref{sec:FisherIID} that the Fisher
information for a distribution over $T$ IID samples is $T$ times the
Fisher information of a single sample.

It can also be shown that the only invariant volume element that can
be defined over a family of distributions is one proportional to the
square root of the determinant of the Fisher information matrix 
$|g|$; {\em i.e.} $\sqrt{|g|}\prod_id\theta_i$.  This can
be regarded as proportional to the number of distinct distributions in
the coordinate prism with opposite corners at $\theta$ and $\theta+d\theta$.
To see this, consider coordinates that diagonalize the Fisher
information matrix in the neighborhood of $\theta$,
so that $dD=(1/2)\sum_ig_{ii}d\theta_i^2$.  Recall
that we consider two distributions distinct if $dD>\alpha/T$, {\em
i.e.} if $(1/2)\sum_ig_{ii}d\theta_i^2>\alpha/T$.  We can therefore
densely pack the distributions by placing them on the corners of a
rectangular lattice separated by $\sqrt{2\alpha/Tg_{ii}}$ along
rectangular coordinate $\theta_i$.  The volume allocated to each
distribution is then $\prod_i\sqrt{2\alpha/Tg_{ii}}$ and the number of
distributions in a rectangular prism with side lengths $d\theta_i$ is
$\left.\prod_id\theta_i\right/\prod_i\sqrt{2\alpha/Tg_{ii}}\propto{\sqrt{|g|}\prod_id\theta_i}$.

\subsubsection{Derivation of Fisher metric from
 Delta Divergence}
\label{sec:FisherFromD}
Consider (\ref{eq:Def:deltaDiv}) with $P(x)=P(x|\theta)$ and
$Q(x)=P(x|\theta+d\theta)$.  Then
\begin{align}
\ln Q(x) &= \ln P(x|\theta+d\theta) \notag\\ 
  &\approx
  \ln P(x|\theta) + \sum_i\partial_{\theta_i}\ln P(x|\theta) d\theta_i
  + \frac{1}{2}\sum_{ij}\partial_{\theta_i}\partial_{\theta_i}\ln P(x|\theta)
  d\theta_i d\theta_j
\end{align}
so abbreviating $P(x|\theta)$ as $P$ and $\partial_{\theta_i}$ as
$\partial_i$ and using
$e^{\epsilon}{\approx}1+\epsilon+\frac{1}{2}\epsilon^2$ gives
\begin{multline}
  P(x)^{\delta}Q(x)^{1-\delta} \approx P^{\delta} e^{({1-\delta})\ln P}
  e^{(1-\delta)\left[\sum_i\partial_i\ln P d\theta_i
  + \frac{1}{2}\sum_{ij}
    \partial_i\partial_j\ln P d\theta_i d\theta_j\right]} \\
  \approx P^{\delta}P^{1-\delta} \left[1+
 (1-\delta) \sum_i\partial_i\ln P d\theta_i
  + \frac{(1-\delta)}{2}
    \sum_{ij}\partial_i\partial_j\ln P d\theta_i d\theta_j \right. \\
   \left. + \frac{(1-\delta)^2}{2}
   \sum_i\partial_i\ln P d\theta_i \sum_j\partial_j\ln P d\theta_j
\right] \\
  \approx P +
 (1-\delta) P \sum_i\partial_i\ln P d\theta_i
  + \frac{(1-\delta)}{2}
    \sum_{ij}\partial_i\partial_j P d\theta_i d\theta_j \\
   - \frac{\delta(1-\delta)}{2}
   P \sum_i\partial_i\ln P d\theta_i \sum_j\partial_j\ln P d\theta_j
\label{eq:Fisherderiv1}
\end{multline}
where we have used
\begin{align}
  \partial_i\partial_j\ln P = 
  \partial_i P^{-1}\partial_j\ln P =
  P^{-1}\partial_i\partial_j P - P^{-2} \partial_i P \partial_j P \notag\\
  = P^{-1}\partial_i\partial_j P - (\partial_i\ln P)(\partial_j\ln P)
\label{eq:d2lnp}
\end{align}
and $(1-\delta)^2-(1-\delta)=(1-\delta)(1-\delta-1)=-\delta(1-\delta)$.
Note that because $\int_xP(x|\theta)=1$, we have 
$\int_x\partial_{\theta_i}P(x|\theta)=0$ and
$\int_x\partial_{\theta_i}\partial_{\theta_j}P(x|\theta)=0$ 
\begin{align}
  \int_xP\partial_i{\ln}P= \int_x\partial_iP=0.
\end{align}
Inserting (\ref{eq:Fisherderiv1}) into (\ref{eq:Def:deltaDiv}) and
applying these identities then gives
\begin{align}
  D_{\delta}[P,Q] &= \frac{1}{\delta(1-\delta)}
  \left[1-\int_x P(X=x)^{\delta} Q(X=x)^{1-\delta}\right] \notag\\
  &\approx \frac{1}{2} \sum_{ij} P(x|\theta)
   \partial_{\theta_i}\ln P(x|\theta)
   \partial_{\theta_j}\ln P(x|\theta) d\theta_i d\theta_j
  = \frac{1}{2} \sum_{ij} g_{ij}d\theta_i d\theta_j
\label{eq:fisherIinfoMatDerived}
\end{align}
in agreement with (\ref{eq:Def:gij}) and (\ref{eq:FisherDist}).

Note that this result is independent of $\delta$.  Dependence on
$\delta$ begins with the 3rd order terms.  These can be expressed by
the {\em Eguchi relations}~\cite{amari93} in
terms of the affine connection coefficients of the $\delta$-geometry.

\subsubsection{Direct Derivation of Fisher information 
for IID data}
\label{sec:FisherIID}
Consider a data set consisting of $T$ IID data points $X=\{x_1,...,x_T\}$.
The likelihood of this data is $P(X|\theta)=\prod_tP(x_t|\theta)$.
The Fisher information is
\begin{equation}
	g_{ij}^{(T)} = \sum_X P(X|\theta) \dthi\ln P(X|\theta)\dthj\ln P(X|\theta)
\end{equation}
Using $\dthi\ln P(X|\theta)=\sum_t\dthi\ln P(X|\theta)$, this is
\begin{equation}
	g_{ij}^{(T)} = \sum_{x_1...x_T} \left[\prod_{t''} P(x_{t''}|\theta)\right]
	          \sum_{tt'} \dthi\ln P(x_t|\theta)\dthj\ln P(x_t'|\theta)
\end{equation}
Observe that factors in the product over $t''$ for which $t''$ is
equal neither to $t$ nor $t'$ are dependent on $x_{t''}$ only 
through an overall factor of $P(x_{t''}|\theta)$, which sums to 1
due to normalization.  The terms with $t\neq t'$ can be factored
into the form
\begin{equation}
	\left[\sum_x P(x|\theta)\dthi \ln P(x|\theta)\right]
	\left[\sum_{x'} P(x'|\theta)\dthj \ln P(x'|\theta)\right] =
	\left[\sum_x \dthi P(x|\theta)\right]
	\left[\sum_{x'} \dthj P(x'|\theta)\right]
\end{equation}
which vanishes due to the normalization condition.  This leaves the
$t=t'$ terms, of which there are $T$.  Hence
\begin{equation}
	g_{ij}^{(T)} = T\sum_x P(x|\theta)\dthi\ln P(x|\theta)\dthj\ln P(x|\theta)
\end{equation}
Hence, we observe that the Fisher information of an IID data set
is simply the data set size times the Fisher information of a single
sample:
\begin{equation}
	g_{ij}^{(T)} = T g_{ij}
\end{equation}

\subsection{Categorical geometry}
The family of categorical distributions over $n$ possible outcomes has $n-1$
independent coordinates $\theta=\{\theta_1,...,\theta_{n-1}\}$,
in terms of which the probability of outcome $x$ is
\begin{equation}
	P(x|\theta) = \theta_x
\end{equation}
where $x$ takes values in $\{1,...,n\}$ and we define
$\theta_n = 1-\sum_{i=1}^{n-1}\theta_i$.  The parameter space is
bounded by the constraints $\theta_i\geq{0}$ and the normalization
constraint $\sum_x\theta_x=1$.  In a subsequent subsection we will
derive the formulas
\begin{align}
	g_{ij} = \frac{\delta_{ij}}{\theta_i} + \frac{1}{\theta_n}
\label{eq:gijMultinomial}
\end{align}
and 
\begin{align}
  g = \prod_{x=1}^n \theta_x^{-1}
\label{eq:gMultinomial}
\end{align}
The invariant volume element is therefore proportional to 
$\sqrt{g}=\prod_x\theta_x^{-1/2}$, which we recognize as the
Dirichlet distribution with all parameters set to $\onehalf$.
Therefore we can determine the normalization constant and determine
that the invariant distribution over categorical distributions is
\begin{align}
  P(\theta|J) =
	\frac{\Gamma(\sum_x\frac{1}{2})}{\prod_x\Gamma(\frac{1}{2})}
				\prod_x\theta_x^{-1/2}
  =	\frac{\Gamma(\frac{n}{2})}{\pi^{\frac{n}{2}}}\prod_x\theta_x^{-1/2}
\label{eq:Jeffreys}
\end{align}
This is also Jeffreys' prior, which motivates the notation $P(\theta|J)$.

\subsection{Derivation of Categorical Fisher Determinant}
This problem can be organized so as to obtain the
characteristic equation as well as the Fisher determinant.  We
introduce eigenvalue parameter $\lambda$, and obtain the Fisher
determinant by setting $\lambda=0$.

In the natural coordinates of the categorical simplex, 
the Fisher information matrix (\ref{eq:gijMultinomial}) has the form
\begin{align}
|g_{\lambda}|=
\begin{vmatrix}
a_1+a_n-\lambda & a_n & a_n & \dots & a_n \\
a_n & a_2+a_n-\lambda & a_n & \dots & a_n \\
a_n & a_n & a_3+a_n-\lambda & \dots & a_n \\
\hdotsfor{5} \\
a_n & a_n & a_n & \dots & a_{n-1}+a_n-\lambda \\
\end{vmatrix} \label{eq:FisherDet1}
\end{align}
where $a_i=\theta_i^{-1}$ for $0{\leq}i{\leq}n$, which includes 
$a_n=\theta_n^{-1}=\left[1-\sum_{i=1}^n\theta_i\right]^{-1}$.
Subtracting the last row of (\ref{eq:FisherDet1}) from each of the
other rows gives
\begin{align}
|g_{\lambda}|=
\begin{vmatrix}
a_1-\lambda & 0 & 0 & \dots & \lambda-a_{n-1} \\
0 & a_2-\lambda & 0 & \dots & \lambda-a_{n-1} \\
0 & 0 & a_3-\lambda & \dots & \lambda-a_{n-1} \\
\hdotsfor{5} \\
a_n & a_n & a_n & \dots & a_{n-1}+a_n-\lambda \\
\end{vmatrix}. \label{eq:FisherDet2}
\end{align}
Next, we add each of columns $1,\dots,n-2$, with column $i$ 
scaled by\footnote{
Here we assume $a_i$ is not an eigenvalue for any $i$.  If this were
so, we would have $\sum_{j}g_{kj}v_j=a_i v_k$ for some
eigenvector $v$.  With $g_{kj}=a_j\delta_{jk}+a_n$ this condition is
$a_kv_k+a_nV=a_iv_k$ with $V=\sum_{j=1}^{n-1}v_j$.  Then
$(a_k-a_i)v_k=a_nV$.  This holds for all $k$ including $k=i$, so with
$\theta$ restricted to the interior of the simplex, we have $V=0$.
But then we have $a_k=a_i$ for all $k$ with nonzero $v_k$, of which
there must be at least 1 (and therefore at least two to enable
$V=0$).  Therefore $a_i$ is not an eigenvalue except perhaps in
special cases where two or more probabilities in the categorical are equal.
}
$\frac{a_{n-1}-\lambda}{a_i-\lambda}$, to column $n-1$.  This produces
\begin{align}
|g_{\lambda}|=
\begin{vmatrix}
a_1-\lambda & 0 & 0 & \dots & 0 \\
0 & a_2-\lambda & 0 & \dots & 0 \\
0 & 0 & a_3-\lambda & \dots & 0 \\
\hdotsfor{5} \\
a_n & a_n & a_n & \dots & a_{n-1}+a_n-\lambda
 + \sum_{i=1}^{n-2}\frac{(a_{n-1}-\lambda)a_n}{a_i-\lambda}
\end{vmatrix}. \label{eq:FisherDet3}
\end{align}
Expanding by minors around the last column, the determinant is now
trivially seen to be the product of the diagonal elements. 
We obtain
\begin{align}
|g_{\lambda}|&=\left[\prod_{i=1}^{n-2}(a_i-\lambda)\right] 
(a_{n-1}-\lambda)\left[1+a_n\sum_{i=1}^{n-1}\frac{1}{a_i-\lambda}\right]
\notag\\
&= \left[\prod_{i=1}^{n-1}(a_i-\lambda)\right] 
\left[1+a_n\sum_{i=1}^{n-1}\frac{1}{a_i-\lambda}\right]
\label{eq:chardet}
\end{align}
in which the $a_n$ term has been brought into the sum as term $i=n-1$.

With $\lambda=0$ we have 
$\sum_{i=1}^{n-1}\frac{1}{a_i-\lambda}=\sum_{i=1}^{n-1}\theta_i=1-\theta_n$
and therefore $1+a_n\sum_{i=1}^{n-1}\frac{1}{a_i}=1+a_n(1-a_n^{-1})=a_n$
so the Fisher determinant is simply
\begin{align}
|g|=|g_0|=\prod_{i=1}^n\theta^{-1}
\end{align}
as was stated in (\ref{eq:gMultinomial}).

\section{Verification of integrals in the main text}
\label{sec:appB}

In this section we verify some of the results given for
integrals appearing in the main text.

\subsection{Integral (\ref{eq:intglA})}
\label{sec:intglA}
To verify (\ref{eq:intglA}),
recall the standard result $d\tan^{-1}(x)/dx=1/(1+x^2)$ to obtain
\begin{align}
  \frac{d}{d\xi} & \biggl[\tan^{-1}(\sqrt{\xi/(1-\xi)})
    -\sqrt{\xi(1-\xi)}\biggr] \notag\\
  &= [1+\xi/(1-\xi)]^{-1}
       \onehalf\sqrt{(1-\xi)/\xi}\frac{(1-\xi)+\xi}{(1-\xi)^2}
       -\onehalf \sqrt{\xi(1-\xi)}^{-1}(1-2\xi) \notag\\
       &= \onehalf (1-\xi) \frac{1}{(1-\xi)^2}\sqrt{\frac{1-\xi}{\xi}}
       -\onehalf (1-2\xi)\frac{1}{\sqrt{\xi(1-\xi)}}  \notag\\
       &= \onehalf \frac{1}{\sqrt{\xi(1-\xi)}}
       -\onehalf (1-2\xi)\frac{1}{\sqrt{\xi(1-\xi)}}
       = \sqrt{\frac{\xi}{1-\xi}} \label{eq:intglAchk}
\end{align}

\subsection{Integral (\ref{eq:indefIntEllip})}
\label{sec:indefIntEllip}
To verify (\ref{eq:indefIntEllip}), let us write the indefinite
integral piecemeal as
\begin{align}
  I &= c K(\sin^{-1}(\xi),m) \\
  c &= 2\frac{b-a}{|b-a|} \frac{1}{\sqrt{b(1-b)}} \\
  \xi &= \sqrt{\frac{(1-b)(b-x)}{(1-a)(a-x)}} \label{eq:xiDef}\\
  m &= \frac{a(1-a)}{b(1-b)} \label{eq:mDef}.
\end{align}
With the substitution
\begin{align}
  t&=\sin\theta \\
  dt &= \cos\theta d\theta \\
  d\theta &= dt/cos\theta = dt/\sqrt{1-\sin^2\theta} = dt/\sqrt{1-t^2},
\end{align}
the definition (\ref{eq:elipDef}) becomes
\begin{align}
  K(\phi,m) &= \int_0^{\sin\phi} dt/\sqrt{(1-t^2)({1-mt^2})}.
\end{align}
so
\begin{align}
  K(\sin^{-1}(\xi),m) &= \int_0^{\xi} dt/\sqrt{(1-t^2)({1-mt^2})}.
\end{align}
Then
\begin{align}
  \frac{dI}{dx} = c\frac{dK(\sin^{-1}(\xi),m)}{d\xi}\frac{d\xi}{dx} =
  c\frac{1}{\sqrt{(1-\xi^2)({1-m\xi^2})}}\frac{d\xi}{dx}. \label{eq:dIdx}
\end{align}
Differentiating (\ref{eq:xiDef}), we have
\begin{align}
  \frac{d\xi}{dx} &= \sqrt{\frac{1-b}{1-a}}\frac{1}{2}
  \sqrt{\frac{a-x}{b-x}} \frac{-(a-x)+(b-x)}{(a-x)^2} \notag\\
  &= \frac{1}{2}\sqrt{\frac{1-b}{1-a}}
  \sqrt{\frac{a-x}{b-x}} \frac{b-a}{(a-x)^2}
\end{align}
Substituting (\ref{eq:xiDef}) into $(1-\xi^2)(1-m\xi^2)$ gives
\begin{align}
  (1-\xi^2)(1-m\xi^2) &= 1-(1+m)\xi^2+m\xi^4 \notag\\
  &= 1-(1+m)\frac{(1-b)(b-x)}{(1-a)(a-x)}
  + m\biggr(\frac{(1-b)(b-x)}{(1-a)(a-x)}\biggl)^2  \notag\\
  &= \frac{\nu}{(1-a)^2(a-x)^2} \\
  \nu &= (1-a)^2(a-x)^2-(1+m)(1-a)(a-x)(1-b)(b-x) \notag\\
  &\qquad + m(1-b)^2(b-x)^2  \label{eq:nuDef}
\end{align}
We observe from (\ref{eq:mDef}) that
\begin{align}
  1+m = \frac{b(1-b)+a(1-a)}{b(1-b)}
\end{align}
so we can write (\ref{eq:nuDef}) as
\begin{align}
  \nu &= (1-a)^2(a-x)^2 -(b(1-b)+a(1-a))(1-a)(a-x)(b-x)/b \notag\\
  &\qquad + a(1-a)(1-b)(b-x)^2/b \notag\\
  b\nu/(1-a) &= (1-a)(a-x)^2b -(b(1-b)+a(1-a))(a-x)(b-x) \notag\\
  &\qquad + a(1-b)(b-x)^2 \notag\\
  &= (1-a)b(a^2-2ax+x^2)+a(1-b)(b^2-2bx+x^2) \notag\\
  &\qquad -(b-b^2+a-a^2))(ab-(a+b)x+x^2) \notag\\
  &= a^2b-a^3b+ab^2-ab^3-ab^2+ab^3-a^2b+a^3b \notag\\
  &\qquad +(-2ab+2a^2b-2ab+2ab^2 \notag\\
  &\qquad\qquad +ab-ab^2+a^2-a^3+b^2-b^3+ab-a^2b)x \notag\\
  &\qquad +(b-ab+a-ab-b+b^2-a+a^2)x^2 \notag\\
  &= (-a^3-b^3 +a^2b+ab^2 +a^2+b^2 -2ab)x
  +(a^2+b^2-2ab)x^2 \notag\\
  &= (-(a^2-b^2)(a-b) + (a-b)^2)x +(a-b)^2x^2 \notag\\
  &= (a-b)^2x^2 + (-(a+b)(a-b)^2 + (a-b)^2)x \notag\\
  &= (a-b)^2x(1+x-a-b).
\end{align}
Then we can write (\ref{eq:dIdx}) as
\begin{align}
  \frac{dI}{dx} &= 2\frac{b-a}{|b-a|} \frac{1}{\sqrt{b(1-b)}}
  \frac{(1-a)(a-x)}{\sqrt{\nu}}
  \frac{1}{2}\sqrt{\frac{1-b}{1-a}}
  \sqrt{\frac{a-x}{b-x}} \frac{b-a}{(a-x)^2} \notag\\ 
  &= \frac{|b-a|\sqrt{1-a}}{\sqrt{b\nu(a-x)(b-x)}}
  = \sqrt{\frac{(b-a)^2}{(a-b)^2x(1+x-a-b)(a-x)(b-x)}} \notag\\
  &= 1/\sqrt{x(a-x)(b-x)(1+x-a-b)}
\end{align}
(remembering that $a-x\geq 0$ and $b-x\geq 0$), in agreement
with the integrand of (\ref{eq:abMrg}).

\section{Some numerical integration details.}
\label{app:Num}

\begin{figure}
  \begin{center}
    \includegraphics[scale=0.6,
      trim=0 0 50 0,clip]{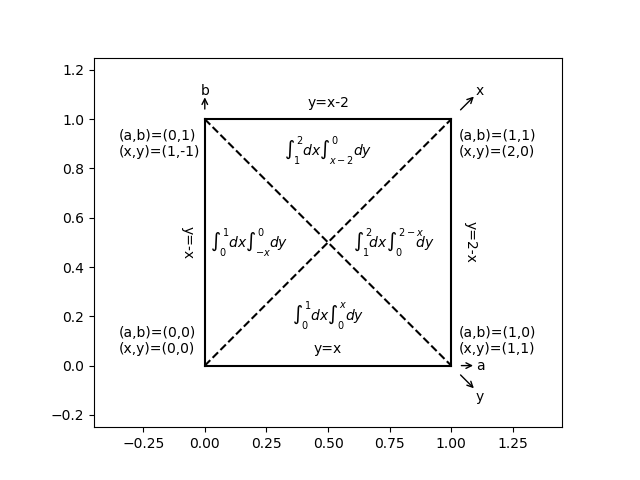}
  \end{center}
  \caption{The singularity-aligned coordinate system
    (\ref{eq:coordTransX}), (\ref{eq:coordTransY}).}
\label{fig:coordTrans}
\end{figure}

For numerical integration
purposes, it is convenient to use a coordinate system in which the
singularities of the prior along the diagonals of the (a,b)
coordinate square lie at constant coordinate values, so that a
convenient variably-sized rectangular grid can be used to sample more
densely near the singularities.  To this end, let
\begin{align}
  x &= a+b &  a=(x+y)/2 \label{eq:coordTransX} \\
  y &= a-b &  b=(x-y)/2 \label{eq:coordTransY} \\
  \frac{\partial(x,y)}{\partial(a,b)} &=
  \begin{vmatrix}
    1 & 1 \\
    1 & -1
  \end{vmatrix}
  = -2 &\\
  a=0 &\iff y=-x & a=1 \iff y=2-x \\
  b=0 &\iff y=x & b=1 \iff y=x-2
\end{align}
so the singularities lie at $x=1$ and $y=0$, and in these coordinates,
the measure must be divided by 2.  We can write
\begin{align}
  2\int_0^1da\int_0^1db\cdots &= \int_0^1dx\int_0^xdy\cdots
  + \int_0^1dx\int_{-x}^0dy\cdots \notag\\
  & + \int_1^2dx\int_0^{2-x}dy\cdots
  + \int_1^2dx\int_{x-2}^0dy\cdots
\end{align}
This coordinate transformation and the regions to which each of these
four integrals apply are illustrated in Figure~\ref{fig:coordTrans}.

The integrals were performed over a non-uniform coordinate-aligned
grid in the $(x,y)$ coordinates, using the same non-linear spacing of
grid points in each dimension, in an attempt to sample more densely
near the singularities.  The interval $[0,2]$ was divided by $2n$ cuts,
with the $i\th$ cut at
\begin{align}
  1-x_{n-i-1}^{\alpha}  & \qquad\qquad 0<i<n \notag \\
  1+x_{i-n}^{\alpha}  & \qquad\qquad n<i<2n \label{eq:sampling} \\
  x_i = \frac{i+0.5}{n} & \qquad\qquad 0\leq i\leq n-1. \notag
\end{align}
These cuts define $2n-1$ adjacent 1D cells, each centered at the
midpoint of its defining cuts.  The function samples were taken at
these midpoints.  For $0<\alpha<1$, this results in more dense sampling
near the ends of the interval than the center.  We chose $\alpha=0.5$,
somewhat arbitrarily.  These intervals were separately
rescaled and recentered to
cover each of the 4 triangles shown in Figure~\ref{fig:coordTrans},
as well as their reflections through each triangle's edge along
boundary of the domain of integration.  Cell centers falling with
the domain of integration were counted; the others dropped.

It should not be difficult to improve upon this sampling strategy,
perhaps replacing it with something alltogether different such as
a Gibbs sampler.  The method used should sample fairly accurately
near the singularities along the diagonal lines, but neglects
to take any special care near the singularities at the external
boundaries.  In addition to not sampling very densely near the
midpoints of those boundaries, no account is taken of the slicing
through the sampling rectangles that cross those boundaries.
Even so, this was adequate for our present purposes.  With
$n=1000$, the integral of the non-informative prior came out to 0.9989,
so the integrals should be accurate to about 0.1\%.

The numerical integration to determine $\kappa(a,b)$
was performed for each $(a,b)$ value by scaling (\ref{eq:sampling})
with $n=100,000$ and $\alpha={0.5}$
to the unit interval and accepting those points lying within
the integration limits of (\ref{eq:Kappadef}).
For the figures, the $(a,b)$ grid was defined with $n=250$.
For Figure~\ref{fig:kappa}, the numerical integration
was performed with $n$ increased to $1,000,000$.

\bibliographystyle{plain}
%\bibliography{\bibd/rr}
\bibliography{}

\end{document}